\documentclass[10pt]{amsart}
\usepackage{amsmath, amssymb, amsthm, a4wide}
\usepackage{geometry}
\usepackage{fancyhdr}
\usepackage{enumerate}
\newtheorem{theorem}{Theorem}[section]
\newtheorem{lemma}[theorem]{Lemma}
\newtheorem{corollary}[theorem]{Corollary}
\newtheorem{definition}[theorem]{Definition}
\newtheorem{proposition}[theorem]{Proposition}
\newtheorem{example}[theorem]{Example}   
\newtheorem{remark}[theorem]{Remark}

\numberwithin{equation}{section}

\newcounter{newlist} 

\newcounter{nnnewlist} 

\newcounter{nelist}

\begin{document}

\title[Lyapunov-type conditions and GSDEs]{Lyapunov-type conditions and stochastic differential equations driven by $G$-Brownian motion}
\author{Xinpeng LI}
\author{Xiangyun LIN}
\author{Yiqing LIN}
\address{Xinpeng LI\newline\indent Institute for advanced research and school of mathematics\newline\indent Shandong University\newline\indent 250100 Jinan, China\newline}
\address{Xiangyun LIN\newline\indent College of Mathematics and Systems Science\newline\indent Shandong University of Science and Technology\newline\indent 266590 Qingdao, China\newline}
\address{Yiqing LIN\newline\indent Fakult\"at f\"ur Mathematik, Universit\"at Wien\newline\indent 1090 Wien, Austria\newline}
\address{\emph{Corresponding author:} {\rm Xiangyun LIN: xiangyun.lin.sdust@gmail.com}}

%\thanks{This research is partially supported by the Post-graduate Study Abroad Program sponsored by China Scholarship Council.}


\date{\today}


\begin{abstract}
This paper studies the solvability and the stability of stochastic differential equations driven by $G$-Brownian motion (GSDEs). In particular, the existence and uniqueness of the solution for locally Lipschitz GSDEs is obtained by localization methods, also the stability of such GSDEs are discussed with Lyapunov-type conditions.
\end{abstract}
\thanks{Xipeng Li research is supported by Fundamental Research Funds of Shandong University under grant No. 2014GN007 and by China Postdoctoral Science Foundation under grant No. 2014M561907; Xiangyun Lin research is supported by the National Natural Science Foundation of China under grant No. 11271007; Yiqing LIN research is supported by European Research Council (ERC) under grant FA506041.}
\maketitle

\section{Introduction}
\noindent Recently, Peng \cite{P3} set up a framework of a time consistent sublinear expectation associated with a new type of Brownian motion that draws a lot of attention. This sublinear expectation, named $G$-expectation, is primarily defined as a sublinear functional $\mathbb{E}[\cdot]$ on a well-posed space $C_{b, lip}(\Omega)$, which is the collection of all bounded and Lipschitz functions 
on the canonical space $\Omega:= \mathcal{C}[0, \infty)$. Under the $G$-expectation, the canonical process $B$ is seen to be a $G$-Brownian motion whose increments are zero-mean, independent, and $G$-normal distributed. Observed in a later work of Denis-Hu-Peng \cite{DHP}, this type of sublinear expectations can be regarded as the supreme of linear expectations over a weakly compact subset of martingale measures and thus, induces a capacity on $\Omega$. Accordingly, the notion of ``quasi-sure'' (q.s.) related to this capacity is brought to this framework. \\

%the the theory is established  and Peng on the bounded and continuous function space ${C}_b(\Omega)$, which is called $G$-expectation.  %that equipped with a capacity induced by the ``upper probability'' over a set of martingale measures. Naturally, under which
%
%``finite dimensional distribution'' of the canonical process can be characterized by a generalized heat equation, and the process with such property is named $G$-Brownian motion.\\
%called $G$-expectation,
As the foundation of the $G$-stochastic calculus, Peng also introduces the $G$-It\^o integral in his pioneer work \cite{P3}. Indeed, this integral is first defined pathwisely as the Riemann-Stieltjes sum on $[0, T]$ for all ``adapted'' step process formed by random variables in $C_{b, lip}(\Omega)$. And then, owing to the $G$-It\^o inequality, this definition can be extended to the $M^2_G$-norm (similar to the classical $M^2$-norm) closure of the step process space. Following Peng's work, related topics in the $G$-stochastic analysis are systematically developed by many authors.\\

In particular, some authors are interested in the stochastic differential equation driven by $G$-Brownian motion (GSDE), which has a similar form as its counterpart in the classical framework, however, holds in a q.s. sense:
\begin{equation}\label{ert}
X_t=x+\int^t_0 f(s, X_s)ds+\int^t_0 h(s, X_s)d\langle B\rangle_s +\int^t_0 g(s, X_s) dB_s,\ 0\leq t\leq T,\ q.s.,
\end{equation}
where $\langle B\rangle$ is the quadratic variation of the $G$-Brownian motion. Under the Lipschitz assumptions on the coefficients, Peng \cite{P3} and Gao \cite{G} have proved the wellposedness of such equation with the fixed-point iteration. Moreover, Bai-Lin \cite{BL} has studied the case when coefficients are integral-Lipschitz and Lin \cite{L1, L2} considers the reflected GSDEs with some good enough boundaries. We observe in these papers the Lipschitz structure of the coefficients and the condition on the reflecting boundary make sure that the solution $X$ and $g(\cdot, X_\cdot)$ live in the space of good integrands defined by Peng \cite{P3}.\\

In contrast to the existing results on GSDEs, the focus of the present paper is on (\ref{ert}) with local Lipschitz coefficients under a Lyapunov-type condition, while detailed discussion to this kind of classical stochastic differential equations (SDEs) can be found, for example, in Has'minski\v{i} \cite{HA}. Following a standard procedure, we shall apply the localization technique to approximate the solution of this equation, however, one can see that the localized coefficient will be eliminated from the integrand space defined in \cite{P3}, since it may be no longer ``continuous'' in $\omega$ somehow. Thanks to the work of Li-Peng \cite{LP}, the space of suitable integrands are essentially expanded, which requires less regularity and only ``local integrability''. Indeed, Li-Peng first extend the definition of the $G$-It\^o integral for all step process formed by bounded random variable, and verify that in their case the $G$-It\^o inequality still holds true, so that this definition can be extended by the completion under the $M^2_G$-norm. Due to this extension, one can consider the $G$-It\^o integral for such process ${\bf 1}_{[0, \tau]}(\cdot)$, where $\tau$ is an $\mathcal{F}^B$-stopping time taking values in $[0, T]$ and thus, the $G$-It\^o integral is said to be well defined in \cite{LP} on an even larger class $M^2_w([0, T])$ by stopping techniques. Taking advantage of these results, we deduce that this kind of GSDEs are well posed on $M^2_w([0, T])$. Particularly, we define the truncated GSDEs that are uniform Lipschitz, and carefully choose the stopping times in order to construct a consistent localized sequence. Finally, the Lyapunov-type condition ensures that the solution can be approximated pathwisely.\\

Another task of this paper is to consider the stability of the GSDEs under perturbations of the initial value, which shows a rough picture of asymptotic behavior of GSDEs under certain conditions. We recall Has'minski\v{i} \cite{HA} and Mao \cite{M} for details in this topic for SDEs. In the $G$-framework, quasi-sure exponential stability of a linear GSDE has been considered in Zhang-Chen \cite{ZC} and the (exponential) moment stability of GSDEs has been studied by Hu-Ren-Xu \cite{HRX} with a particular Lyapunov function (independent of $h$ and $g$). In what follows, we obtain the sufficient conditions for the moment stability by means of a general Lyapunov function and investigate the relation between the moment stability and the quasi-sure stability for GSDEs. Moreover, a property of $G$-Brownian motion is obtained at the end of this paper.\\

This paper is organized as follows: Section 2 provides preliminaries in the $G$-framework for our further discussion; in Section 3 we conclude the main idea in \cite{LP} and explain the stopping techniques in the $G$-framework; Section 4 studies the solvability of local Lipschitz GSDEs while Section 5 shows the stability results for such equations; in the appendix, Proposition 1.4 in §IV-1 of [12] is generalized, which is essential for our argument.
\section{Preliminaries}\label{32s}
We first recall some notions and results of $G$-expectation and the related spaces of random variables and stochastic processes. The reader interested
in more detailed description of these notions and results is referred to
Denis-Hu-Peng \cite{DHP}, Gao \cite{G}, Li-Peng \cite{LP} and Peng \cite{P3}.\\[6pt]

Let $\Omega$ be the space of all $\mathbb{R}^{d}$-valued
continuous paths with $\omega_0 =0$
equipped
with the distance
$$
\rho(\omega^1, \omega^2):=\sum^\infty_{N=1} 2^{-N} ((\max_
{t\in[0,N]} | \omega^1_t-\omega^2_t|)  \wedge 1),
$$
and let $B_t(\omega):=\omega_t$ be the canonical process.\\

Set
$$
L_{ip}(\Omega):=\{\varphi(B_{t_1}, \cdots, B_{t_n}):n\geq
1,\ 0\leq t_1\leq \cdots\leq t_n<\infty,\ \varphi \in C_{b,lip}
(\mathbb{R}^{d\times n})\},
$$
where $C_{b,lip}
(\mathbb{R}^{d\times n})$ the collection of all bounded and Lipschitz functions on $\mathbb{R}^{d\times n}$.\\

Fixing a sublinear monotone function $G$: $\mathbb{S}^d\rightarrow \mathbb{R}$, the related $G$-expectation on $(\Omega, L_{ip}(\Omega))$ can be constructed in the following way: for each $\xi\in L_{ip}(\Omega)$ with the form
$$\xi=\varphi(B_{t_1},B_{t_2}-B_{t_1},\cdots,B_{t_n}-B_{t_{n-1}}), \ \ 0\leq t_1<t_2<\cdots<t_n,$$
define
$$\mathbb{E}[\xi]:=u_1(0,0),$$
where $u_1(0,0)\in\mathbb{R}$ is obtained via the following procedure: for $k=n,\cdots,1$, $u_k:=u_k(t,x;x_1,\cdots,x_{k-1})$ is a function of $(t,x)$ with the parameter $(x_1,\cdots,x_{k-1})\in\mathbb{R}^{d\times(k-1)}$, which is the solution of the following $G$-heat equation defined on $[t_{k-1},t_k)\times\mathbb{R}^d$:
 $$\frac{\partial u_k}{\partial t}-G(D^2 u_k)=0$$
 with terminal condition
 $$u_k(t_k,x;x_1,\cdots,x_{k-1})=u_{{k+1}}(t_{k},x;x_1,\cdots,x_{k-1},x),$$
where $u_{n+1}(t_n,x;x_1,\cdots,x_{n-1},x):=\varphi(x_1,\cdots,x_{n-1},x)$. Note that under this $\mathbb{E}[\cdot]$, the canonical process $(B_t)_{t\geq 0}$ is a $G$-Brownian motion (see \S III-2 in \cite{P3}).\\

We denote by $L_G^p(\Omega)$ the completion of $Lip(\Omega)$ under the natural norm $||\cdot||_p:=\mathbb{E}[|\cdot|^p]^{\frac{1}{p}}$. For each $t\in \lbrack0,\infty)$, we list the following notations:

\begin{itemize}
\item[$\bullet$] $\Omega_{t}:=\{ \omega_{\cdot \wedge t}:\omega \in \Omega \}$;
\item[$\bullet$] $\mathcal{F}_{t}:=\mathcal{B}(\Omega_{t})$;
\item[$\bullet$] $L^{0}(\Omega)$: the space of all $\mathcal{B}(\Omega)$-measurable real functions;
\item[$\bullet$] $L^{0}(\Omega_{t})$: the space of all $\mathcal{B}(\Omega_{t}%
)$-measurable real functions;
\item[$\bullet$] $B_{b}(\Omega)$: all bounded elements in $L^{0}(\Omega)$; $B_{b}%
(\Omega_{t}):=B_{b}(\Omega)\cap L^{0}(\Omega_{t})$;
\item[$\bullet$] $C_{b}(\Omega)$: all continuous elements in $B_{b}(\Omega)$;
$C_{b}(\Omega_{t}):=C_{b}(\Omega)\cap L^{0}(\Omega_{t})$.\\
\end{itemize}

We can extend the domain of $G$-expectation $\mathbb{E}[\cdot]$ from $L_{ip}(\Omega)$ to $L^{0}(\Omega)$ by the procedure introduced in \cite{DHP},
i.e., constructing an upper expectation $\bar{\mathbb{E}}[\cdot]$:
\begin{equation*}
\bar{\mathbb{E}}[X]:=\sup_{\mathbb{P}\in\mathcal{P}_G}E^\mathbb{P}[X],\ X\in L^0(\Omega),
\end{equation*}
where $\mathcal{P}_G$ is a weakly compact family of martingale measures on $(\Omega,\mathcal{B}(\Omega))$. This upper expectation
coincides with the $G$-expectation $\mathbb{E}[\cdot]$ on $L_{ip}(\Omega)$ and thus, on its completion $L^1_G(\Omega)$. Naturally, a Choquet capacity can be defined by:
$$
\bar{C}(A):=\sup_{\mathbb{P}\in\mathcal{P}_G}\mathbb{P}(A),\ A\in\mathcal{B}(\Omega),
$$
and the notation of ``quasi-surely'' (q.s.) is introduced as follows:
\begin{definition}
A set $A\in\mathcal{B}(\Omega)$ is called polar if $\bar{C}(A)=0.$ A property
is said to hold quasi-surely if it holds outside a polar
set.
%We denote by $\mathcal{N}^{\bar{C}}(\Omega_T)$ all the polar sets in $\mathcal{B}(\Omega_T)$ and we set $\bar{\mathcal{B}}(\Omega_t):=\mathcal{B}(\Omega_t)\vee \mathcal{N}^{\bar{C}}(\Omhega_T)$ \footnote{\noindent If we append the following class of sets into $\mathcal{B}(\Omega_t)$, the argument throughout this chapter will not alter: $\mathcal{N}^{\bar{C}}:=\{A\subset B: B\ {\rm is\ a}\ \mathbb{P}{\rm -null\ set\ in\ \mathcal{B}(\Omega_T)},\ {\rm for\ all}\ \mathbb{P}\in\mathcal{P}_G\}.
%$
%}.
\end{definition}

We note that the following Borel-Cantelli Lemma, Markov's inequality and Upwards convergence theorem also holds for such Choquet capacity and upper expectation  (Lemma 5, Lemma 13 and Theorem 10 in \cite{DHP})
\begin{lemma}\label{MI}
Let $X\in L^0(\Omega)$ satisfying $\bar{\mathbb{E}}[|X|^p]<\infty$, for $p>0$. Then, for each $a>0$,
$$
\bar{C}(\{|X|>a\})\leq \frac{\bar{\mathbb{E}}[|X|^p]}{a^p}.
$$
\end{lemma}
\begin{lemma}\label{BCL}
Let $\{A_n\}_{n\in \mathbb{N}}\subset \mathcal{B}(\Omega)$ such that
$$
\sum^{\infty}_{n=0}\bar{C}(A_n)<\infty.
$$
Then, $\limsup_{n\rightarrow \infty}A_n$ is a polar set.
\end{lemma}
\begin{theorem}\label{CTU}
Let $\{X^n\}_{n\in\mathbb{N}}\subset L^0(\Omega_T)$ be a sequence such that $X^n\uparrow X$, q.s., and there exists a $\mathbb{P}\in\mathcal{P}_G$ with $E^\mathbb{P}[X^0]>-\infty$, then $\bar{\mathbb{E}}[X^n]\uparrow\bar{\mathbb{E}}[X]$.
\end{theorem}

We also have Fatou's lemma (Lemma 2.11 in Bai-Lin \cite{BL}) in the $G$-framework:
\begin{lemma}\label{fatulem}
Assume that $\{X^n\}_{n\in\mathbb{N}}$ is a sequence in $L^0(\Omega)$ and for a $Y\in L^0(\Omega)$ satisfying $\bar{\mathbb{E}}[|Y|]<\infty$ and all $n\in\mathbb{N}$, $X^n \geq Y$, q.s., then
$$
\bar{\mathbb{E}}[\liminf_{n\rightarrow\infty}X^n]\leq\liminf_{n\rightarrow\infty}\bar{\mathbb{E}}[X^n].\\
$$
\end{lemma}

In order to introduce stopping time in $G$-framework, Li-Peng \cite{LP} first consider the process space generalized by $B_b(\Omega)$ and define the related It\^{o}'s integral. Let $L^p_*(\Omega)$ be the completion of $B_b(\Omega)$, and consider the
following type of simple processes:
\begin{align*}
M^0_{b}([0,T])  &  =\{ \eta:\eta_{t}(\omega)=\sum_{i=0}^{N-1}\xi_{i}%
(\omega)\mathbf{1}_{[t_{i},t_{i+1})}(t),\\
&  \forall N\in \mathbb{N},\ 0=t_{0}<\cdots<t_{N}=T,\ \xi_{i}\in B_{b}%
(\Omega_{t_{i}}),\ i=0,\cdots,N-1\}.
\end{align*}

\begin{definition}
For $p\geq 1$, we denote by $M_*^p([0,T])$ the completion of $M_b^{0}([0,T])$ under the following norm:
\begin{equation*}
\|\eta\|_p:=\bigg(\bar{\mathbb{E}}\bigg[\frac{1}{T}\int_0^T|\eta_t|^p dt\bigg]\bigg)^{1/p}.
\end{equation*}
\end{definition}
%\indent Proved in \cite{P3}, the $G$-It\^o type integral can be first defined on $M_c^{0}([0,T])$ as a Riemann sum which is contracting so that can be extended to $M_G^p([0,T])\rightarrow L^p_G(\Omega_T)$. Afterwards, Li-Peng \cite{LP} observe that $G$-It\^o's inequality (Proposition 2.11 in \cite{LP}) still holds for all elements in $M_b^0([0,T])$, then the methodology mentioned above can be applied for defining the $G$-It\^o type integral on $M_*^p([0,T])$.\\[6pt]
%\noindent

Here below is the definition of the $G$-It\^o integral. In the sequel, $B^{\bf a}$ denotes the inner product of ${\bf a}\in \mathbb{R}^d$ and $B$, which is still a $G$-Brownian motion, and $\sigma_{\bf aa^{\rm T}}:=\bar{\mathbb{E}}[({\bf a}, B_1)^2]$.

\begin{definition}\label{hao}
%Defining $:=({\bf a}, B)$, for a fixed, $B^{\bf a}$  is a scalar $G$-Brownian motion.
%Giving , we note $.
%be a given vector in $$ and , where $({\bf a}, B)$ denotes the scalar product of ${\bf a}$. Following Li and Peng \cite{LP},
For each $\eta\in M^0_b([0, T])$, we define the It\^o type integral
$$
\mathcal{I}_{[0, T]}(\eta)=\int^T_0\eta_tdB^{\bf a}_t:=\sum^{N-1}_{k=0}\xi_k(B^{\bf a}_{t_{k+1}}-B^{\bf a}_{t_k}).
$$
Then, the linear mapping $\mathcal{I}_{[0, T]}$ on $M^0_b([0, T])$ can be continuously extended to $\mathcal{I}_{[0, T]}: M^2_*([0, T])\rightarrow L^2_*(\Omega_T)$ and for each $\eta \in M^2_*([0, T])$, we define $\int^T_0\eta_tdB^{\bf a}_t:=\mathcal{I}_{[0, T]}(\eta)$.
\end{definition}

We notice that Gao \cite{G} proves the $G$-It\^{o} integral $X_\cdot=\int_0^\cdot\eta_sdB^{\bf a}_s$ has a continuous $\bar{C}$-modification, for any $\eta\in M_G^2(0,T)\subset M_*^2([0,T])$, where $M_G^2([0,T])$ is defined in \cite{P3} by the completion of simple process formed by random variables in $C_{b, lip}(\Omega)$. In fact, this result still holds for $\eta\in M_*^2(0,T)$, that is to say, there exists $\bar{X}$ which is continuous in $t$ pathwisely and $\bar{C}(|X_t-\bar{X}_t|\neq 0)=0$ for all $t\in [0,T]$. Furthermore, we can deduce the following proposition, so that we do not care the choice of the $t$-continuous $\bar{C}$-modification for the $G$-It\^o integral $(\int^t_0 \eta_sdB_s)_{0\leq t\leq T}$.
\begin{proposition}\label{indis}
Suppose that $X$ and $X'$ are two processes that have $t$-continuous paths and for each $t\in [0, T]$, $X_t=X'_t$, q.s.. Then, $X$ and $X'$ are indistinguishable in the q.s. (indistinguishable, for short in the rest part of this paper).
\end{proposition}

\section{Localisation techniques and extension of the $G$-It\^o type integral}\label{s33}\noindent
% and Their $G$-Stochastic Integrals}
\noindent In this section, we follow Li-Peng \cite{LP} to define the $G$-It\^o type integral for a certain type of locally integrable processes and generalize some results in \cite{LP}. The notion of stopping times always refer to the ones with respect to $\mathcal{F}$.
%Notice that we have enlarged a little the filtration, however, this will not alter any result in \cite{LP}.
% For some local integrable processes, $G$-stochastic integrals have also been defined by some stopping time technique in .
%\noindent We now give the definition to an admissible class of  processes that are local integrable.
\begin{definition}\label{MPW}
For fixed $p\geq 1$, a stochastic process $\eta$ is said to be in $M^p_w([0, T])$, if it is associated with a sequence $\{\sigma_m\}_{m\in\mathbb{N}}$ of stopping times incresing to infinity such that
$$\eta{\bf 1}_{[0, \sigma_m\wedge T]} \in M^p_*([0, T]),\ \ \forall m\in\mathbb{N}.$$
\end{definition}
%\begin{remark}
%Given $\tau$ that is another stopping time, for each $m\in\mathbb{N}$, it is easy to verify that $\eta {\bf 1}_{[0, \sigma_m\wedge\tau]}\in M^p_*([0, T])$.
%\end{remark}
\begin{remark}
One can easily deduce from this definition that
\begin{equation}\label{haha}
\int^T_0|\eta_t|^pdt<\infty,\ q.s..
\end{equation}
\end{remark}
This definition is a slight modification of original one in \cite{LP}. We can define the related It\^{o}'s integral in the similar way.\\

\indent For a given $\eta\in M^2_w([0, T])$ associated with $\{\sigma_m\}_{m\in\mathbb{N}}$, we note $\tau_m:=\sigma_m\wedge T$ and consider
%Lemma 4.2 in Li and Peng \cite{LP}, for
%each $m\in\mathbb{N}$,
the $t$-continuous $\bar{C}$-modification of
$(\int^t_0 \eta_s{\bf 1}_{[0, \tau_m]}(s)dB^{\bf a}_s)_{0\leq t\leq T}$.
% is well defined process that has $t$-continuous paths.
For each $m$, $n\in\mathbb{N}$, $n >m$, by Lemma 4.2 in \cite{LP} and the continuity of the paths of the $G$-It\^o integral, we can find a polar set $\hat{A}^{m, n}$, such that for all $\omega\in (\hat{A}^{m, n})^c$, the following equality holds:
\begin{align}\label{keypoint}
\int^{t\wedge\tau_m}_0\eta_sdB^{\bf a}_s(\omega)
&=\int^{t}_0\eta_s{\bf 1}_{[0, \tau_m]}(s)dB^{\bf a}_s(\omega)\notag\\
&=\int^{t}_0\eta_s{\bf 1}_{[0, \tau_m]}(s){\bf 1}_{[0, \tau_n]}(s)dB^{\bf a}_s(\omega)\\
&=\int^{t\wedge\tau_m}_0\eta_s{\bf 1}_{[0, \tau_n]}(s)dB_s^{\bf a}(\omega),\ 0\leq t\leq T.\notag
\end{align}
Define a polar set
\begin{equation*}
\hat{A}:=\bigcup_{m=1}^{\infty}\bigcup_{n=m+1}^{\infty}\hat{A}^{m, n}.
\end{equation*}
For each $m\in\mathbb{N}$ and $(\omega, t)\in\Omega\times[0, T]$, we set
$$
X^m_t(\omega):=
\left\{
\begin{aligned}
\int^{t}_0\eta_s{\bf 1}_{[0, \tau_m]}(s)dB^{\bf a}_s(\omega)\ &,\
\omega\in \hat{A}^c\cap \bar{\Omega};\\
0\ &,\ {\rm otherwise}.
\end{aligned}\right.
$$
For each $\omega\in \hat{A}^c$ and $m$, $n\in\mathbb{N}$, $n>m$, $X^n(\omega)\equiv X^m(\omega)$ on $[0, \tau_m(\omega)]$. Therefore, we can define unambiguously a process by stipulating that it equal to $X^m$ on $[0, \tau_m(\omega)]$.

\begin{definition}\label{locdef} Giving $\eta\in M^2_w([0, T])$,
for each $(\omega, t)\in\Omega\times[0, T]$, we define
\begin{equation*}
\int^{t}_0\eta_sdB^{\bf a}_s(\omega):=\lim_{m\rightarrow\infty}X^m_t(\omega).
\end{equation*}
\end{definition}
\begin{remark}
For $\eta\in M^p_w([0, T])$, where $p\geq 2$, associated with $\{\sigma_m\}_{m\in\mathbb{N}}$,
%We would like to give an example of the process in $M^2_w([0, T])$.
define $X_t:=\int^t_0\eta_sdB_s$, then $X$ stays in $M^p_w([0, T])$, where the sequence of stopping times can be chosen as $\{\sigma_m\}_{m\in\mathbb{N}}$.
%we prove that $\sigma_m:=\inf\{t\geq 0: |X_t|\geq m\}$ is a suitable sequence such that with respect to this sequence, $X$ satisfies the two conditions in Definition \ref{MPW}. First, $X$ has continuous paths outside a polar set $A$, thus for each $\omega\in A^c$, there exists an $m\in\mathbb{N}$ such that $\sigma_m(\omega)=T$, which implies . Second,
\end{remark}
\indent At the end of this section, we generalize the BDG type inequality in \cite{LP} for $\eta\in M^2_w([0, T])$.
\begin{lemma}\label{le5}
Let
%\footnote{After proving that $\{\int^t_0 \eta_s dB_s\}_{0\leq t\leq T}$
%has a $t$-continuous $\bar{C}$-modification, one can follow Lin \cite{L2} to generalize this result to the case that $0<p<\infty$, by directly applying the BDG type inequality under each $\mathbb{P}\in\mathcal{P}_G$. However, the continuity of paths of $\int\eta dB$ has not been mentioned before the proof of the BDG type inequality in Lin \cite{L2}.} either,
$\eta\in M^p_w([0, T])$, the following inequality holds: for $p>0$, there exists $C_p>0$ depends only on $p$, such that
\begin{equation*}
\bar{\mathbb{E}}\bigg[\sup_{0\leq s\leq t}\bigg|\int^s_0\eta_udB^{\bf a}_u\bigg|^p\bigg]\leq C_p\sigma^{p/2}_{\bf aa^{\rm T}}\bar{\mathbb{E}}\bigg[\bigg(\int^t_0|\eta_s|^2ds\bigg)^{p/2}\bigg].
\end{equation*}
%independent of ${\bf a}$, $\eta$ and $\Gamma$.
\end{lemma}
\noindent {\bf Proof:} Suppose $\eta$ is associated with stopping times $\{\sigma_m\}_{m\in \mathbb{N}}$, then we have
\begin{equation*}
\bar{\mathbb{E}}\bigg[\sup_{0\leq s\leq t}\bigg|\int^s_0\eta_u{\bf 1}_{[0, \sigma_m\wedge T]}dB^{\bf a}_u\bigg|^p\bigg]\leq C_p\sigma^{p/2}_{\bf aa^{\rm T}}\bar{\mathbb{E}}\bigg[\bigg(\int^t_0|\eta_s|^2{\bf 1}_{[0, \sigma_m\wedge T]}ds\bigg)^{p/2}\bigg]\leq C_p\sigma^{p/2}_{\bf aa^{\rm T}}\bar{\mathbb{E}}\bigg[\bigg(\int^t_0|\eta_s|^2ds\bigg)^{p/2}\bigg].
\end{equation*}
Since $\int^\cdot_0\eta_tdB^{\bf a}_t$ has continuous paths,
$$
\sup_{0\leq s\leq t}\bigg|\int^s_0\eta_u{\bf 1}_{[0, \sigma_m\wedge T]}dB^{\bf a}_u\bigg|^p=\sup_{0\leq s\leq t}\bigg|\int^{s\wedge \sigma_m\wedge T}_0\eta_u dB^{\bf a}_u\bigg|^p\rightarrow \sup_{0\leq s\leq t}\bigg|\int^{s}_0\eta_udB^{\bf a}_u\bigg|^p,\ q.s..
$$
Finally, Fatou's lemma (Lemma \ref{fatulem}) gives desired result.\hfill$\square$
\section{Locally Lipschtiz $G$-stochastic differential equations}\label{s35}
\noindent In this section, we study the solvability of the following $n$-dimensional stochastic differential equation driven by $d$-dimensional $G$-Brownian motion: 
%on $M^2_*([0, T]; \mathbb{R})$ 
%of the following form:
%We present our main results in this section. First, we redo the job of Gao \cite{G} and Peng \cite{P3} to show that the solvability of GSDE still holds, although 
%Then, we proved the existence and uniqueness results for multidimensional RGSDEs in a bounded convex domain and in some more general domains.
% by giving the existence and uniqueness of the solutions to the scalar valued RGSDEs with Lipschitz coefficients. Additionally, a comparison theorem is given at the end of this paper.
%\subsection{$G$-Stochastic differential equations}
%\noindent We consider the following $n$-dimensional GSDE:
\begin{equation}\label{gsde1}
X_t=x+\int^t_0 f(s, X_s)ds+\int^t_0 h(s, X_s)d\langle B, B\rangle_s +\int^t_0 g(s, X_s) dB_s,\ 0\leq t\leq T,\ q.s.,
\end{equation}
where $x\in\mathbb{R}^n$ is the initial value and $\langle B, B\rangle=(\langle B^i, B^j\rangle)_{i, j = 1, \ldots, d}$ is the mutual variation matrix of $B$.\\[6pt]
% and $f$, $h$ and $g$ are functions such that for fixed $\omega\in\Omega$, 
%$f(\cdot, \cdot)(\omega): [0, T]\times\mathbb{R}\rightarrow \mathbb{R}$,  
%$h(\cdot, \cdot)(\omega)=(h^{ij}(\cdot, \cdot)(\omega))_{i,j=1, \ldots, d}:[0, T]\times\mathbb{R}\rightarrow \mathbb{R}^{d^2}$ and $g(\cdot, \cdot)(\omega)=(g^{j}(\cdot, \cdot)(\omega))_{j=1, \ldots, d}:[0, T]\times\mathbb{R}\rightarrow \mathbb{R}^{d}$.\\[6pt]
%We emphasize that all results in this section can be generalized to the multi-dimensional case with no difficulty.\\[6pt]
\indent In a first stage, we consider this equation under the so-called uniformly Lipschitz conditions:\\[6pt]
%tion $x\in \mathbb{R}$;
{\bf (H1)} For some $p\geq 2$ and each $x\in\mathbb{R}^n$, $f(\cdot, x)$, $h^{ij}(\cdot, x)$, $g^{j}(\cdot, x)\in M^p_*([0,T]; \mathbb{R}^n)$, $i, j=1, \ldots, d$;\\[3pt]
{\bf (H2)} The coefficients $f$, $h$ and $g$ are uniformly Lipschitz in $x$, i.e., for each $t\in[0, T]$ and $x, x'\in\mathbb{R}^n$, 
\begin{align*}|f(t,x)-f(t, x')|+||h(t, x)-h(t, x')||
+||g(t, x)-g(t, x')||\leq C_L|x-x'|,
\end{align*}
where $||\cdot||$ is the Hilbert-Schmidt norm of a matrix. \\[6pt]
%Here, q.s. means this inequality holds for all the $\omega$ outside a polar set $A$ independent of $t$.\\[3pt]
\indent We notice that the domain of coefficients here is a little larger than the ones in Gao \cite{G} and Peng \cite{P3}, however, the same method as in \cite{G} can be applied to prove the uniqueness and existence theorem as follows:
\begin{theorem}\label{tgsde}
Let (H1) and (H2) hold. Then,  for any $p\geq 2$, there exists a unique process $X\in M^p_*([0, T]; \mathbb{R}^n)$ that has $t$-continuous paths and satisfies the GSDE (\ref{gsde1}). Denote by $X^x$ and $X^y$ the solutions starting with $x$, $y\in\mathbb{R}^n$, then there exists $C>0$ that depends only on $p$, $T$ and $C_L$, such that
\begin{equation*}
\bar{E}[\sup_{t\in[0, T]}|X^x_t-X^y_t|^p]\leq C|x-y|^p.
\end{equation*}
\end{theorem}
\begin{remark}
Similarly to Lemma 5.1 in Bai-Lin \cite{BL}, we have the following assertion: suppose that for some $p\geq 1$, $\zeta$ is a function such that  
$\zeta(\cdot, x)\in M^p_*([0,T]; \mathbb{R}^n)$, for each $x\in \mathbb{R}^n$ and that $\zeta(\cdot, \cdot)$ is uniformly Lipschitz in $x$,
%condition, i.e., for each $t\in[0, T]$ and each $x_1$, 
%$x_2\in \mathbb{R}^n$, $|\zeta(t, x_1)-\zeta(t, x_2)|\leq C_L|x_1-x_2|$.
then for each $X\in M^p_*([0,T]; \mathbb{R}^n)$, $\zeta(\cdot, X_\cdot)$ is an element in $M^p_*([0,T]; \mathbb{R}^n)$. This result ensures that (\ref{gsde1}) is well defined. 
%We shall prove the existence part by the Picard iteration. 
%Before proceeding this, we should prove the following lemma, which ensures that all stochastic integrals are well defined in each iterative step. In the sequel, a constant $C>0$ that depends only on $p$, $n$, $T$ and $C_L$ may vary from line to line.
\end{remark}
\indent In what follows, we consider such a GSDE, whose coefficients satisfy both a locally Lipschitz condition and a Lyapunov's type condition. Here below are two assumptions that we concern: from now on, we adopt the Einstein notation. \\[6pt]
{\bf (H2')} The coefficients $f(\cdot, \cdot)$, $h^{ij}(\cdot, \cdot)$, $g^{j}(\cdot, \cdot): [0, T]\times\mathbb{R}^n\rightarrow \mathbb{R}^n$ are deterministic functions continuous in $t$ and locally Lipschitz in $x$, i.e., for each $x$, $x'\in\{a: |a|\leq {R}\}$, there exists a positive constant $C_R$ that depends only on $R$, such that for each $t\in [0, T]$,
\begin{align*}|f(t, x)-f(t, x')|+||h(t, x)-h(t, x')||
+||g(t, x)-g(t, x')||\leq C_R|x-x'|.
\end{align*}
{\bf (H3')} There exist a deterministic nonnegative Lyapunov function $V\in \mathcal{C}^ {1, 2}([0, T]\times\mathbb{R}^n)$, such that 
$$\inf_{|x|\geq R}\inf_{t\in[0, T]}V(t, x)\rightarrow \infty,\ {\rm as}\ R\rightarrow \infty,$$
and for some constant $C_{LY}>0$ and all $(t, x)\in[0, T]\times \mathbb{R}^n$,
$$
\mathcal{L}V(t, x)\leq C_{LY}V(t, x),
$$ 
where $\mathcal{L}$ is a differential operator defined by
\begin{align*}
%\label{oper}
\mathcal{L}V
%&:= \bigg(\partial_{x^\nu} V(x) f^{\nu}(x)+\sup_{\mathcal{S}\in \Sigma}
%\big(\partial_{x^\nu} V(x) h^{\nu}_{ij}(x)\sigma^\mathcal{S}_{ij}+\frac{1}{2}\partial_{x^\mu x^\nu}^2 V(x) g^{\mu}_{i}(x)g^{\nu}_{j}(x)\sigma^\mathcal{S}_{ij}\big)\bigg)\\
=\partial_{t} V+\partial_{x^\nu} V f^{\nu}+G\bigg(\left(\partial_{x^\nu} V\cdot (h^{\nu ij}+h^{\nu ji})+\partial_{x^\mu x^\nu}^2 V\cdot g^{\mu i}g^{\nu j}\right)^d_{i, j=1}\bigg).
\end{align*}
%in which $\mathcal{S}:=(\sigma^{\mathcal{S}}_{ij})_{i, j =1}^d\in\Sigma\subset \mathbb{S}^d$.
\begin{remark}
The coefficients $f$, $h^{ij}$, $g^j$, $i, j=1, \ldots, d$, are continuous in $t$ and thus uniformly continuous on $[0, T]$, so for each $x\in \mathbb{R}^n$, $f(\cdot, x)$, $h^{ij}(\cdot, x)$, $g^j(\cdot, x)\in M^p_G([0, T]; \mathbb{R}^n)\subset M^p_*([0, T]; \mathbb{R}^n) $, for any $p\geq 1$. 
\end{remark}
\begin{lemma}\label{lm5}
For each $X\in M^1_w([0, T]; \mathbb{R}^n)$ that has continuous paths and for each function $V\in \mathcal{C}^ {1, 2}([0, T]\times\mathbb{R}^n)$,  we have $\eta (V, X)\in M^p_w([0, T]; \mathbb{S}^d)$, for any $p\geq 1$, where $$\eta^{ij}_\cdot(V, X):= \partial_{x^\nu} V(\cdot, X_\cdot) (h^{\nu ij}(\cdot, X_\cdot)+h^{\nu ji}(\cdot, X_\cdot))+\partial_{x^\mu x^\nu}^2 V(\cdot, X_\cdot) g^{\mu i}(\cdot, X_\cdot)g^{\nu j}(\cdot, X_\cdot).$$
\end{lemma}
\noindent {\bf Proof:} Similarly to the proof of Theorem 5.4 in Li-Peng \cite{LP}, we can find a sequence of stopping times $\{\sigma_m\}_{m\in\mathbb{N}}$ satisfying Definition \ref{MPW}, such that on $[0, \sigma_m]$, $|X|$ is bounded by $m$. It is obvious that $X{\bf 1}_{[0, \sigma_m]}\in M^p_*([0, T]; \mathbb{R}^n)$, for any $p\geq 1$.\\[6pt]
\indent Noticing that $\partial_{x^\mu}V$, $\partial^2_{x^\mu x^\nu}V$, $h^{ij}$ and $g^j$, $\mu$, $\nu=1,\ldots n$,  $i, j=1,\ldots, d$, are bounded and uniformly continuous on compact sets, we can approximate these functions uniformly on $[0, T]\times B(0, m)$ by $C_{b, lip}([0, T]; \mathbb{R}^n)$ ones and thus, we have $\eta^{ij}_{\cdot\wedge \sigma_m}(V, X)\in M^p_*([0, T])$, for any $p\geq 1$, which implies the desired result. \hfill$\square$
\begin{theorem}\label{lf}
Let (H2') and (H3') hold. Then, for any $p\geq 2$, there exists a unique solution of the GSDE (\ref{gsde1}) in $M^p_w([0, T]; \mathbb{R}^n)$ that has $t$-continuous paths on $[0, T]$, and the following estimate holds:
$$
\bar{\mathbb{E}}[V(t, X^x_t)]\leq e^{C_{LY}T}V(0, x).
$$
%For two $x$, $y\in\mathbb{R}^n$, let $X^x$ and $X^y$ are two solution of (\ref{gsde1}) with the initial value $x$ and $y$, then there exist a constant $C>0$ that depends only on $p$, $n$, $T$ and $C_L$ such that
%\begin{equation}\label{ju}
%\bar{E}[\sup_{t\in[0, T]}|X^x_t-X^y_t|^p]\leq C|x-y|^p.
%\end{equation}
\end{theorem}
\noindent\textbf{Proof:} {\it Existence.}
For each $N\in\mathbb{N}$, we first consider the following truncated GSDE:
\begin{equation}\label{eqtrun}
X^N_t=x+\int^t_0 f^N(s, X^N_s)ds+\int^t_0 h^N(s, X^N_s)d\langle B, B\rangle_s +\int^t_0 g^N(s, X^N_s) dB_s,\ 0\leq t\leq T,\ q.s.,
\end{equation}
where $f^N$, $(h^{ij})^N$ and $(g^{j})^N$, $i$, $j=1, \ldots, d$, are defined in the following way:
\begin{align*}
\zeta^N(t, x)=
\left\{\begin{array}{c@{,}l}
\zeta(t, x)&\ {\rm{if}}\ |x|\leq N;\\[3pt]
\zeta({t, Nx}/{|x|})&\ {\rm{if}}\ |x|> N.
\end{array}\right.
\end{align*}
It is easy to verify that $f^N$, $h^N$ and $g^N$ are all bounded functions and uniformly Lipschitz in $x$. Then, by the result of Lipschitz GSDEs with coefficients in $M^p_G([0, T]; \mathbb{R}^n)$, for example in \cite{G} and \cite{P1},
the truncated GSDE (\ref{eqtrun}) admits a unique solution $X^N\in M^p_G([0, T]; \mathbb{R}^n)$, for any $p\geq 2$, whose paths are $t$-continuous.\\[6pt]
\indent Define a sequence of stopping times by
\begin{equation*}
\tau_N:=\inf\{t: |X^N_t|\geq N\}\wedge T,
\end{equation*}
which satisfies $\{\tau_N\leq t\}\in \mathcal{F}_t$.
Thanks to Lemma \ref{ere}, we can deduce from (\ref{eqtrun}) that
\begin{align*}
X^N_{t\wedge \tau_N}&=x+\int^t_0 f^N(s, X^N_s){\bf 1}_{[0, \tau_N]}(s)ds\\
&+\int^t_0 h^N(s, X^N_s){\bf 1}_{[0, \tau_N]}(s)d\langle B, B\rangle_s +\int^t_0 g^N(s, X^N_s){\bf 1}_{[0, \tau_N]}(s) dB_s\\
&=x+\int^t_0 f^{N+1}(s, X^N_s){\bf 1}_{[0, \tau_N]}(s)ds\\
&+\int^t_0 h^{N+1}(s, X^N_s){\bf 1}_{[0, \tau_N]}(s)d\langle B, B\rangle_s +\int^t_0 g^{N+1}(s, X^N_s){\bf 1}_{[0, \tau_N]}(s) dB_s\\
&=x+\int^{t\wedge \tau_N}_0 f^{N+1}(s, X^N_s)ds\\
&+\int^{t\wedge \tau_N}_0 h^{N+1}(s, X^N_s)d\langle B, B\rangle_s+\int^{t\wedge \tau_N}_0 g^{N+1}(s, X^N_s)dB_s,\
0\leq t\leq T,\ q.s..
\end{align*} 
On the other hand, by the definition of $X^{N+1}$, we have
\begin{align}\label{trad}
X^{N+1}_{t\wedge \tau_N}=x+\int^{t\wedge \tau_N}_0 f^{N+1}(s, X^{N+1}_s)ds&+\int^{t\wedge \tau_N}_0 h^{N+1}(s, X^{N+1}_s)d\langle B, B\rangle_s\notag\\ &+\int^{t\wedge \tau_N}_0 g^{N+1}(s, X^{N+1}_s) dB_s,\ 0\leq t\leq T,\ q.s..
\end{align}
By the uniqueness of the solution to the truncated GSDE (\ref{trad}), for each $N\in \mathbb{N}$,
%with coefficients $f^{N+1}$, $h^{N+1}$ and $g^{N+1}$, %whose components are all in $M^p_*([0, T])$ for fixed $x$ and uniformly Lipschitz in $x$,
$X^N$ and $X^{N+1}$ are distinguishable on $[0, \tau_N]$. This also implies that the sequence $\{\tau_N\}_{N\in\mathbb{N}}$ are q.s. increasing.\\[6pt]
\indent Now we aim to show that 
\begin{equation}\label{ver1}
\bar{C}\bigg(\bigcup_{N=1}^{\infty}
\{\omega: \tau_N(\omega)=T\}\bigg)=1.
\end{equation}
Because $|X^N_{\cdot}|$ never exceeds $N$ on $[0, \tau_N]$, we have
\begin{align}\label{sucheng}
f(t, X^N_{t}){\bf 1}_{[0, \tau_N]}(t)=f^N(t, X^N_t){\bf 1}_{[0, \tau_N]}(t);\ &h(t, X^N_{t}){\bf 1}_{[0, \tau_N]}(t)=h^N(t, X^N_t){\bf 1}_{[0, \tau_N]}(t);\notag\\
&g(t, X^N_{t}){\bf 1}_{[0, \tau_N]}(t)=g^N(t, X^N_{t}){\bf 1}_{[0, \tau_N]}(t),\ 0\leq t\leq T,
\end{align}
where the right-hand side are $M^p_*([0, T]; \mathbb{R}^n)$ processes (cf. Lemma 4.2 in \cite{LP}), for any $p\geq 2$.
As a result, 
\begin{align*}
X^{N}_{t\wedge \tau_N}
&=x+\int^{t\wedge \tau_N}_0 f^N(s, X^N_s)ds+\int^{t\wedge \tau_N}_0 h^N(s, X^N_s)d\langle B, B\rangle_s +\int^{t\wedge \tau_N}_0 g^N(s, X^N_s) dB_s\notag\\
&=x+\int^{t\wedge \tau_N}_0 f^N(s, X^N_{s}){\bf 1}_{[0, \tau_N]}(s)ds\\
&+\int^{t\wedge \tau_N}_0 h^N(s, X^N_{s}){\bf 1}_{[0, \tau_N]}(s)d\langle B, B\rangle_s+\int^{t\wedge \tau_N}_0 g^N(s, X^N_{s}){\bf 1}_{[0, \tau_N]}(s) dB_s\\
&=x+\int^{t}_0 f(s, X^N_s){\bf 1}_{[0, \tau_N]}(s)ds\\
&+\int^{t}_0 h(s, X^N_s){\bf 1}_{[0, \tau_N]}(s)d\langle B, B\rangle_s
+\int^{t}_0 g(s, X^N_s){\bf 1}_{[0, \tau_N]}(s) dB_s,\ 0\leq t\leq T,\ q.s..\notag
\end{align*}
Then, applying $G$-It\^o's formula (cf. Theorem 5.4 in \cite{LP}) to 
$\Phi(t\wedge \tau_N, X^N_{t\wedge\tau_N}):= \exp(-C_{LY}(t\wedge \tau_N))V(t\wedge\tau_N, X^N_{t\wedge\tau_N})$, we have
\begin{align}\label{po}
\Phi(t\wedge \tau_N, X^N_{t\wedge \tau_N})&-\Phi(0, x)\notag\\[3pt]
&=\int^{t\wedge \tau_N}_0(\partial_{t}\Phi(s, X^N_s)+\partial_{x^\nu}\Phi(s, X^N_s)f^{\nu}_{s}(s, X^N_s))ds\notag\\
&+\int^{t\wedge \tau_N}_0\partial_{x^\nu}\Phi(s, X^N_s)g^{\nu j}_{s}(s, X^N_s)dB^j_s\\
&+\int^{t\wedge \tau_N}_0\partial_{x^\nu}\Phi(s, X^N_s)h^{\nu ij}_{s}(s, X^N_s)\notag\\
&+\frac{1}{2}\partial^2_{x^\mu x^\nu}\Phi(s, X^N_s)g^{\mu i}_s(s, X^N_s) g^{\nu j}_s(s, X^N_s)d\langle B^i, B^j\rangle_s,\ q.s..\notag
\end{align}
Borrowing the notation in Lemma \ref{lm5}, the above equality can be written into 
\begin{align*}
Phi(t\wedge \tau_N, X^N_{t\wedge \tau_N})&-\Phi(0, x)\\[3pt]
&=\int^{t\wedge \tau_N}_0\partial_{t}\Phi(s, X^N_s)+\partial_{x^\nu}\Phi(s, X^N_s)f^{\nu}_{s}(s, X^N_s)+G\big(\eta_s(\Phi, X^N)\big)ds\\
&+\int^{t\wedge \tau_N}_0\partial_{x^\nu}\Phi(s, X^N_s)g^{\nu j}_{s}(s, X^N_s)dB^j_s\notag\\
&+\int^{t\wedge \tau_N}_0 \eta^{ij}_s(\Phi, X^N) d\langle B^i, B^j\rangle_s-\int^{t\wedge \tau_N}_0G(\eta_s(\Phi, X^N))ds\notag\\
&=\int^{t}_0\mathcal{L}\Phi(s, X^N_s){\bf 1}_{[0, \tau_N]}(s)ds+\int^{t}_0\partial_{x^\nu}\Phi(s, X^N_s)g^{\nu j}_{s}(s, X^N_s){\bf 1}_{[0, \tau_N]}(s)dB^j_s\notag\\
&+\int^{t}_0 \eta^{ij}_s(\Phi, X^N){\bf 1}_{[0, \tau_N]}(s) d\langle B^i, B^j\rangle_s-\int^{t}_0G(\eta_s(\Phi, X^N){\bf 1}_{[0, \tau_N]}(s))ds,\ q.s..\notag
\end{align*} 
%$X^N_{\cdot}\equiv X^N_{\cdot\wedge \tau_N}$ on $[0, \tau_N]$, we have for each $t\in[0, T]$, $\mathcal{L}\Phi(t\wedge \tau_N, X^N_{t\wedge\tau_N})\leq 0$. 
From (\ref{sucheng}) and the fact that $\partial_x V(t, x)$ is uniformly continuous in $t$ and uniformly Lipschitz in $x$ on $[0, T]\times B(0, N)$, it is readily observed that 
$\partial_\nu \Phi(\cdot, X^N_{\cdot})g^{\nu j}(X^N_\cdot){\bf 1}_{[0, \tau_N]}(\cdot)\in M^p_*([0, T])$, for any $p\geq 2$ (cf. Proposition 3.11 § in \cite{LP}). Then, we obtain
$$
\bar{\mathbb{E}}\bigg[\int^{t}_0\partial_{x^\nu}\Phi(s, X^N_s)g^{\nu j}_{s}(s, X^N_s){\bf 1}_{[0, \tau_N]}(s)dB^j_s\bigg]=0.
$$ 
On the other hand, $\eta (\Phi, X)\in M^p_w([0, T]; \mathbb{S}^d)$, for any $p\geq 1$, which implies $\eta (\Phi, X){\bf 1}_{[0, \tau_N]}\in M^p_w([0, T]; \mathbb{S}^d)$. Here, we claim that 
$$
\bar{\mathbb{E}}\bigg[\int^{t}_0 \eta^{ij}_s(\Phi, X^N){\bf 1}_{[0, \tau_N]}(s) d\langle B^i, B^j\rangle_s
-\int^{t}_0G\big(\eta_s(\Phi, X^N){\bf 1}_{[0, \tau_N]}(s)\big)ds
\bigg]\leq 0,
$$
whose proof is postponed to the appendix (see Lemma \ref{lms}). Hence,
$$
\bar{\mathbb{E}}[\Phi(T\wedge \tau_N, X^{N}_{T\wedge \tau_N})]-\Phi(0, x)\leq\bar{\mathbb{E}}\bigg[\int^{T\wedge \tau_N}_0 \mathcal{L}\Phi(t, X^N_{t})dt\bigg].
$$
Because $\mathcal{L}V\leq C_{LY}V$ implies $\mathcal{L}\Phi\leq 0$, 
\begin{equation*}
\bar{\mathbb{E}}[V(T\wedge \tau_N, X^{N}_{T\wedge \tau_N})]\leq V(0, x)\exp(C_{LY}T).
\end{equation*}
In particular, we have
\begin{equation*}
\bar{\mathbb{E}}[{\bf 1}_{\{\tau_N<T\}}V(T\wedge \tau_N, X^{N}_{T\wedge \tau_N})]\leq V(0, x)\exp(C_{LY}T).
\end{equation*}
Since $X^N$ has $t$-continuous paths, $\tau_N<T$ implies $|X^{N}_{T\wedge \tau_N}|=N$, q.s., from which we deduce
$$
\bar{C}(\{\omega: \tau_N(\omega)<T\}) \inf_{|x|\geq N}\inf_{t\in[0, T]}V(t, x)\leq V(0, x)\exp(C_{LY}T).
$$
As $N\rightarrow \infty$, by (H3'),  we obtain
$$
1\geq \lim_{N\rightarrow \infty}\bar{C}(\{\omega: \tau_N(\omega)=T\})
\geq 1- \lim_{N\rightarrow \infty}\bar{C}(\{\omega: \tau_N(\omega)<T\})=1.
$$
\indent Since $\{\omega: \tau_N(\omega)=T\}$ is increasing, the upwards convergence theorem
%\sup_{N\in\mathbb{N}}\bar{C}(\{\omega: \tau_N(\omega)=T\})=\sup_{N\in\mathbb{N}}\sup_{\mathbb{P}\in\mathcal{P}_G}
%\mathbb{P}(\{\omega: \tau_N(\omega)<T\})
%=\sup_{\mathbb{P}\in\mathcal{P}_G}\sup_{N\in\mathbb{N}}
%\mathbb{P}(\{\omega: \tau_N(\omega)<T\})=1,
%$$
%which 
yields (\ref{ver1}). Therefore, there exists a polar set $A$, such that for all $\omega\in A^c$, the following assertion holds: one can find an $N_0(\omega)$ that depends on $\omega$, such that for all $N\geq N_0(\omega)$, $N\in\mathbb{N}$, $\tau_N(\omega)=T$. Then, we define
\begin{equation}\label{soludef}
X_t(\omega)
=\left\{
\begin{aligned}
X^{N_0(\omega)}_t(\omega), \ 0\leq t\leq T\ &,\
\omega\in A^c;\\
0\ &,\ \omega\in A.
\end{aligned}\right.
\end{equation}
From the argument above, we have for each $\tau_N$, $X{\bf 1}_{[0, \tau_N]}=X^N {\bf 1}_{[0, \tau_N]}\in M^p_*([0, T]; \mathbb{R}^n)$ and thus, $X\in M^p_w([0, T]; \mathbb{R}^n)$, for any $p\geq 2$. Moreover, 
\begin{align*}
X_{t\wedge \tau_N}=X^N_{t\wedge \tau_N}
&=x+\int^{t\wedge \tau_N}_0 f^N(s, X^N_s)ds+\int^{t\wedge \tau_N}_0 h^N(s, X^N_s)d\langle B, B\rangle_s +\int^{t\wedge \tau_N}_0 g^N(s, X^N_s) dB_s\notag\\
&=x+\int^{t\wedge \tau_N}_0 f(s, X_s)ds+\int^{t\wedge \tau_N}_0 h(s, X_s)d\langle B, B\rangle_s\\
&+\int^{t\wedge \tau_N}_0 g(s, X_s) dB_s,\ 0\leq t\leq T,\ q.s.,\notag
%&=x+\int^{t}_0 f(X_s){\bf 1}_{[0, \tau_N]}(s)ds+\int^{t}_0 h(X^N_s){\bf 1}_{[0, \tau_N]}(s)d\langle B, B\rangle_s +\int^{t}_0 g(X^N_s){\bf 1}_{[0, \tau_N]}(s) dB_s,\ 0\leq t\leq T,\ q.s..\notag
\end{align*}
which implies that $X$ satisfies (\ref{gsde1}).\\[6pt]\indent Thanks to the positivity of $C_{LY}$, we have for a fixed $t\in[0, T]$,
\begin{equation*}
\bar{\mathbb{E}}[V(t\wedge \tau_N, X^{N}_{t\wedge \tau_N})]\leq V(0, x)\exp(C_{LY}T).
\end{equation*}
Letting $N\rightarrow \infty$, (\ref{ver1}) and (\ref{soludef}) yield $V(t\wedge \tau_N, X^{N}_{t\wedge \tau_N})
\rightarrow V(t, X_{t})$, q.s., then by Fatou's Lemma we can deduce
\begin{equation*}
\bar{\mathbb{E}}[V(t, X_{t})]\leq V(0, x)\exp(C_{LY}T).
\end{equation*}\\[3pt]
{\it Uniqueness.} Suppose that $X$ and $X'\in M^p_w([0, T]; \mathbb{R}^n)$ are two solutions of the GSDE (\ref{gsde1}) associated with $\{\mu_m\}_{m\in\mathbb{N}}$, and $\{\mu'_m\}_{m\in\mathbb{N}}$. Then, by the continuity of solutions, we can define a new sequence of stopping times by $\tau_m:=\sigma_m\wedge \sigma'_m\wedge\mu_m\wedge\mu'_m$, where 
$$
\sigma_m:=\inf\{t: |X_t|\geq m\}\wedge T,\ \mbox{and}\  \sigma'_m:=\inf\{t: |X'_t|\geq m\}\wedge T,\ m\in \mathbb{N}.
$$ 
From the pathwise uniqueness of Lipschitz GSDEs, we know that $X$ and $X'$ are indistinguishable on $[0, \tau_m]$, for each $m\in \mathbb{N}$. Thus, the uniqueness of the GSDE (\ref{gsde1}) can be deduced from the fact $\bar{C}(\lim_{m\rightarrow \infty}\{\tau_m=T\})=1$.
% if there is another $X'$ 
\hfill{}$\square$
\begin{remark}
The equality (\ref{po}) can be obtained in the following way: first applying $G$-It\^o's formula to $\Phi(t, X^N_t)$, one can see that  (\ref{po}) holds q.s. for each $t\in[0, T]$. Since the processes on both sides of (\ref{po}) have continuous paths, they are indistinguishable. Thus, (\ref{po}) holds q.s. for any bounded stopping time, such as $t \wedge \tau$.
%
%we can work with a bounded the by the continuity of the paths on both sides,During the proof, we have used an extended notion of $G$-It\^o's formula in the following form:
%where $\tau$ is a stopping time. In fact, if $\tau$ equals to the deterministic time $T$, the two sides of (\ref{po}) equal q.s. as a result of Theorem 5.4 in Li and Peng \cite{LP}. %Furthermore, both sides of (\ref{po}) have $t$-continuous paths outside a polar set, then they are . Therefore, for any bounded stopping time $t \wedge \tau$,  (\ref{po}) holds with no problem.
\end{remark}
\begin{example}
In particular, locally Lipschitz GSDEs with a linear growth condition can be regarded as a special case satisfying (H2') and (H3').  For example, letting $n=d=1$, suppose all the coefficients are continuous in $t$, locally Lipschtz in $x$ and for some $K>0$,
$$xf(t, x)+xh(t, x)+ |g(t, x)|^2\leq K(1+|x|^2).$$
In this case, $V$ can be defined by $V(x)=C_K(1+|x|^2)$, where $C_K$ depends on $\Sigma$ and $K$. 
\end{example}
\begin{example}
In fact, (H2') allows us to consider some GSDEs with polynomial growth coefficients. Here is an example, which is similar to Duffing and van der Pol oscillators in random mechanics (see Arnold \cite{AR} for more examples): letting $n=2$ and $d=1$, consider the following GSDE: 
$$
d\left(
\begin{array}{c}
X_t \\
Y_t\\
\end{array}
\right)=\left(
\begin{array}{c}
Y_t\\
-\alpha X_t-\beta X^3_t-\gamma Y_t\\
\end{array}
\right)d\langle B\rangle_t+\left(
\begin{array}{c}
0\\
\sigma\\
\end{array}
\right)dB_t,
$$
where $\alpha$, $\beta$, $\gamma$ and $\sigma$ are positive constants. In this case, the Lyapunov function could be 
$$
V(x, y)=1+\frac{1}{2}y^2+\frac{\alpha}{2}x^2+\frac{\beta}{4}x^4.
$$
\end{example}
\section{Stability of $G$-stochastic differential equations}\label{s36}
\noindent In this section, we study the moment stability of the G-stochastic system (\ref{gsde1}) under perturbations of the initial value. Moreover, we consider the relationship between the moment stability and the q.s. stability under the linear growth condition and obtain a property of the sample paths of the $G$-Brownian motion.\\[6pt]
\indent In order to study the behavior of the solution to the GSDE (\ref{gsde1}) over an infinite horizon, we first give the following definition:
\begin{definition}
We call $X$ a global solution of the GSDE (\ref{gsde1}) if $X$ has $t$-continuous paths on $[0, \infty)$ and for each $T>0$, $X_{\cdot\wedge T}\in M^2_w([0, T]; \mathbb{R}^n)$ is a solution of the GSDE (\ref{gsde1}) on $[0, T]$.
\end{definition}
\begin{remark}
Under (H2') and (H3'), such a global solution exists, since we can first consider the solution of the GSDE (\ref{gsde1}) on $[0, N]$, for each $N\in\mathbb{N}$, and then concatenate all these solutions. The uniqueness can be deduced in the light of the continuity of paths.
\end{remark}
\indent Throughout this section, we always consider the GSDE (\ref{gsde1}) with the locally Lipschitz condition (H2'). Meanwhile, we shall limit ourselves to conditions for stability of the trivial solution $X\equiv 0$. Accordingly, we assume that \\[6pt]
\noindent {\bf (H4')}
$$
f(t, 0)=h(t, 0)=g(t, 0)\equiv 0.
$$
\begin{definition}
Denote by $X^{s, x}$ the solution of the GSDE (\ref{gsde1}) starting with $X_s=x$, $x\in \mathbb{R}^n$. The trivial solution of this $G$-stochastic system in $\mathbb{R}^n$ is said to be
\begin{enumerate}[(i)]
\item $p$-stable, for some $p>0$, if for each $\varepsilon>0$, there exists a $\delta>0$, such that
$$
\sup_{|x|\leq\delta}\sup_{t\geq s}\bar{\mathbb{E}}[|X^{s, x}_t|^p]<\varepsilon;
$$
\item asymptotically $p$-stable, if it is $p$-stable and moreover $$
\bar{\mathbb{E}}[|X^{s, x}_t|^p]\rightarrow 0, {\rm as}\ t\rightarrow \infty;
$$
\item exponentially $p$-stable, if for some positive constants $C$ and $\lambda$
$$
\bar{\mathbb{E}}[|X^{s, x}_t|^p]\leq C|x|^p e^{-\lambda (t-s)}.
$$
\end{enumerate}
In particular, when $p=2$, we call this system is (asymptotically, exponentially) stable in mean square.
\end{definition}
\begin{example}\label{pg}
Letting $n=1$ and $d=1$, we consider the ``population growth model'' in the $G$-framework, i.e. the following GSDE:
\begin{equation}\label{leqe}
X^x_t=x+\int^t_0 \alpha X_s ds+\int^t_0 \beta X_s d\langle B\rangle_s +\int^t_0 \gamma X_s dB_s,
\end{equation}
where $\alpha$, $\beta$ and $\gamma$ are constants. The explicit solution of the above equation is given by Theorem 1.2 in \S V-1 of \cite{P3}:
$$
X^x_t=x\exp\bigg(\alpha t+\big(\beta-\frac{1}{2}\gamma^2\big)\langle B\rangle_t+\gamma B_t\bigg).$$
Thus, for any $p>0$,
$$
|X^x_t|^p= |x|^p \exp\bigg( \alpha p t+\frac{1}{2}p\big(2\beta+\gamma^2(p-1)\big)\langle B\rangle_t\bigg)\mathcal{E}(\gamma pB).
$$
It is obviously that $\gamma p B$ satisfies Assumption 2.1 in Xu-Shang-Zhang \cite{XSZ}, which implies that $\mathcal{E}(\gamma p B):=\exp(\gamma p B-\frac{1}{2}(\gamma p)^2\langle B\rangle)$ is a $G$-martingale. Therefore,
(a) if $2\beta+\gamma^2(p-1)\geq 0$ and $\alpha +\frac{1}{2}\overline{\sigma}^2\big(2\beta+\gamma^2(p-1)\big)<0$, or (b) if $2\beta+\gamma^2(p-1)<0$ and $\alpha +\frac{1}{2}\underline{\sigma}^2\big(2\beta+\gamma^2(p-1)\big)<0$, then the trivial solution of the GSDE (\ref{pg}) is exponentially $p$-stable.
\end{example}
\indent Similar to Theorem \ref{lf}, the sufficient condition for the stability of the $G$-stochastic system (\ref{gsde1}) will be given in terms of Lyapunov functions. Without loss of generality, we consider only the case $s=0$.
\begin{theorem}
Consider the $G$-stochastic system (\ref{gsde1}) with (H2') and (H4'). Suppose that there exists a function $V(t, x)\in\mathcal{C}^{1, 2}([0, \infty)\times \mathbb{R}^n)$ such that for all $t\geq 0$, some positive constant $c_1$ and $c_2$, and some $p>0$,
\begin{equation}\label{54k}
c_1 |x|^p\leq V(t, x)\leq c_2|x|^p.
\end{equation}
Then,
\begin{enumerate}[(a)]
\item  the trivial solution is $p$-stable, if
\begin{equation}\label{54k1}
\mathcal{L}V\leq 0,
\end{equation}
%of the $G$-stochastic system (\ref{gsde1})
\item the trivial solution is exponentially $p$-stable, if there exists a $\lambda>0$ such that for all $(t, x)\in ([0, \infty)\times \mathbb{R}^n)$,
\begin{equation}\label{3}
\mathcal{L}V(t, x)\leq -\lambda V(t,x ). 
\end{equation}
\end{enumerate}
\end{theorem}
\noindent{\bf Proof:} Obviously, (H3') is satisfied in both cases, so that fixing $t>0$, we can follow the procedure in the  proof of Theorem \ref{lf} to construct the solution of  the GSDE (\ref{gsde1}) on $[0, T]$ with the initial value $x$, for some $T>t$. 
%(a) Defining the stopping time $\tau_N$ by (\ref{taun}) and applying $G$-It\^o's formula, we have
%\begin{align}\label{fane}
%V(t\wedge \tau_N, X^x_{t\wedge \tau_N})-V(0, x)
%&=\int^{t\wedge \tau_N}_0\mathcal{L}V(s, X^x_s)ds+\int^{t\wedge \tau_N}_0\partial_{x^\nu}V(s, X^x_s)g^{\nu j}_{s}(s, X^x_s)dB^j_s\notag\\
%&+\int^{t\wedge \tau_N}_0 \eta^{ij}_s(V, X^x) d\langle B^i, B^j\rangle_s-\int^{t\wedge \tau_N}_0G(\eta_s(V, X^x))ds
%,\ q.s..
%\end{align}
%Taking $\bar{\mathbb{E}}[\cdot]$ on both sides of (\ref{fane}), we obtain
%\begin{equation*}\label{keye}
%\bar{\mathbb{E}}[V(t\wedge \tau_N, X^x_{t\wedge \tau_N})-V(0, x_0)]\leq \bar{\mathbb{E}}\bigg[\int^{t\wedge \tau_N}_0 \mathcal{L}V(t, X^x_{t})dt\bigg].
%\end{equation*}
%From (\ref{54k}) and (\ref{54k1}), we deduce
%$$
%\bar{\mathbb{E}}[V(t\wedge \tau_N, X^x_{t\wedge \tau_N})]\leq V(0, x)\leq c_2 |x|^p,
%$$
%then Fatou's Lemma (Lemma \ref{fatulem}) gives 
%$$
%\bar{\mathbb{E}}[V(t, X^x_{t})]\leq V(0, x)\leq c_2 |x|^p.
%$$
%Recalling again (\ref{54k}), we have 
%$$
%\bar{\mathbb{E}}[|X^x_t|^p]\leq \frac{c_2}{c_1} |x|^p,
%$$
%from which one can see that $X\equiv 0$ is $p$-stable.\\[6pt]
%(b) 
Meanwhile, the following inequality can be deduced:
\begin{equation*}
\bar{\mathbb{E}}[V(t, X^x_{t})]\leq V(0, x)\exp(-\lambda t).
\end{equation*}
Finanlly, (\ref{54k}) implies
$$
\bar{\mathbb{E}}[|X^x_t|^p]\leq \frac{c_2}{c_1} |x|^p\exp(-\lambda t),
$$
from which we can deduce both (a) and (b).
\hfill$\square$
\begin{corollary}
Let $p>0$. Consider the stochastic system (\ref{gsde1}) with (H2') and (H4'). Assume that there exists a $P\in \mathbb{S}^n_+$, and constants $\alpha_1\in \mathbb{R}$, $0\leq \alpha_2<\alpha_3$, such that for all $(t, x)\in [0, \infty) \times \mathbb{R}^n$, 
$$
x^{\rm T}Pf(x, t)+G\big((x^{\rm T}P(h^{ij}(x, t)+h^{ji}(x, t)))^d_{i, j=1}+g^{\rm T}(x, t)P g(x, t)\big)\leq \alpha_1 x^{\rm T} P x,
$$
and 
$$
\alpha^2_2(x^{\rm T}Px)^2\leq 2G\big(g^{\rm T}(x, t) P^{\rm T} x x^{\rm T}Pg(x, t)\big)\leq \alpha^2_3(x^{\rm T}P x)^2.
$$
Then, the trivial solution is exponentially $p$-stable if  
(a) $\alpha_1<0$ and $p<2+|\alpha_1|/\alpha_3^2$; or if (b) $0\leq \alpha_1<\alpha^2_2$ and $p<2-2\alpha_1/\alpha_2^2.$
\end{corollary}
\noindent{\bf Proof:} Set $V(x, t):= (x^{\rm T}Px)^{\frac{p}{2}}$. The fact $P\in \mathbb{S}^n_+$ implies
$$
\lambda_{min}^{\frac{p}{2}}(P)|x|^p\leq V(x)\leq \lambda_{max}^{\frac{p}{2}}(P)|x|^p.
$$
By the definition of the differential operator $\mathcal{L}$ and the sub-additivity of the function $G$, we have
\begin{align*}
\mathcal{L}V(x, t)&\leq p (x^{\rm T}Px)^{\frac{p}{2}-1}\bigg(x^{\rm T}Pf(x,t)
+G\big((x^{\rm T}P(h^{ij}(x, t)+h^{ji}(x, t)))^d_{i, j=1}+g^{\rm T}(x, t)P g(x, t)\big)\bigg)\\
&+2p\bigg(\frac{p}{2}-1\bigg)(x^{\rm T}Px)^{\frac{p}{2}-2}G\big(g^{\rm T}(x, t) P^{\rm T} x x^{\rm T}Pg(x, t)\big).
\end{align*}
Then, if (a) holds true (without loss of generality, assuming $p\geq 2$),
$$
\mathcal{L}V(x, t)\leq -p\bigg(|\alpha_1|-\big(\frac{p}{2}-1\big)\alpha^2_3\bigg)V(x, t);
$$
if (b) holds true,
$$
\mathcal{L}V(x, t)\leq -p\bigg(\big(\frac{p}{2}-1\big)\alpha^2_2-\alpha_1\bigg)V(x, t).
$$
In both cases, the exponentially $p$-stability follows from the previous theorem. \hfill$\square$\\[6pt]
%%for some $c_1$, $c_2>0$, for all $x\in \mathbb{R}^n$.
\noindent Now we have the following examples to show the algebraic criteria for the exponential stability in mean square of a linear $G$-stochastic system:
\begin{example}
Consider the following linear GSDE driven by a scalar $G$-Brownian motion:
\begin{equation}\label{ex}
X_t=x+\int^t_0 FX_sds+\int^t_0 HX_s d\langle B\rangle_s+\int^t_0 CX_s dB_s,
\end{equation}
where $F$, $H$ and $C\in \mathbb{R}^{n\times n}$. If there exists a matrix $P\in \mathbb{S}^n_+$, such that for any $x\in \mathbb{R}^n$, 
\begin{equation}\label{ri}
x^{\rm T}\big(PF+ I_n\big)x + G\big(x^{\rm T}(2PH+C^{\rm T}PC)x\big) \leq 0. 
\end{equation}
Then, one can verify that the Lyapunov function $V(x):=x^{\rm T}Px$ satisfies (\ref{54k}) and (\ref{54k1}) and thus, the trivial solution  is exponential stable in mean square.\\[6pt]
\indent Indeed, 
$$
D_xV(x)=2Px,\ D^2_{xx}V(x)\equiv 2P,
$$
thus,
$$
\mathcal{L}V(x)=2\langle Px, Fx\rangle + G\big(\langle 2Px, 2Hx\rangle + \langle 2PCx, Cx\rangle \big).
$$
Obviously, if $P$ satisfies (\ref{ri}), then
$$
\mathcal{L}V(x)\leq =-2|x|^2.
$$
%Moreover,
%$P\in \mathbb{S}^n_+
%$ implies
%$$
%c_1|x|^2\leq V(x)=x^{\rm T} P x\leq c_2|x|^2,
%$$
%for some $c_1$, $c_2>0$, for all $x\in \mathbb{R}^n$. \\[6pt]
\indent In fact, (\ref{ri}) is the Riccati type inequality related to the stability problem of (\ref{ex}) in the $G$-framework, but it is no longer a linear matrix inequality (LMI) because of the non-linearity of the function $G$. However,  a sufficient condition to ensure that $P\in \mathbb{S}^d_+$ satisfies (\ref{ri}) can be given by the following system of LMIs: for some $\alpha\in\mathbb{R}$,
\begin{align*}
\left\{\begin{array}{l}
2PH+C^{\rm T}PC\leq \alpha I_n;\\[3pt]
2PF+I_n\leq -G(\alpha) I_n.
\end{array}\right.
\end{align*}
\end{example}
\indent In what follows, we shall consider a GSDE whose coefficients satisfy a linear growth condition. Similar to the result in the classical framework, the exponentially $p$-stability of the trivial solution to such a GSDE implies the q.s. exponential stability. 
\begin{theorem} \label{mows}
Consider the $G$-stochastic system (\ref{gsde1}) with (H2'), (H4') and the following conditions: for some $K>0$ and all $(t, x)\in [0, \infty)\times \mathbb{R}^n$,
\begin{equation}\label{cone}
x^{\rm T} f(t, x)\vee ||x^{\rm T} h(t,x)||\vee ||g(t,x) ||^2 \leq K|x|^2.
\end{equation}
\noindent If the trivial solution is exponential $p$-stable, for some $p >0$, namely, for all $x\in \mathbb{R}^n$,
\begin{equation}\label{stacond}
\bar{\mathbb{E}}[|X^{s, x}_t|^p]\leq |x|^p \exp(-\lambda (t-s)),
\end{equation}
where $\lambda>0$. 
Then, 
$$
\limsup_{t\rightarrow \infty}\frac{1}{t}\log(|X^{s, x}_t|)\leq-\frac{\lambda}{p}<0,\ q.s.,$$
which implies that the trivial solution is q.s. exponentially stable.
\end{theorem}
\noindent{\bf Proof:} Without loss of generality, we consider only the case $s=0$. Applying $G$-It\^o's formula to $(|X^{x}_t|^2+\delta)^{\frac{p}{2}}$, we have for each $m\in \mathbb{N}$,
\begin{align*}
(|X^x_t|^2+\delta)^{\frac{p}{2}}=(|X^x_{m-1}|^2+\delta)^{\frac{p}{2}}&+\int^t_{m-1} p(|X^x_u|^2+\delta)^\frac{p-2}{2}(X^x_u)^{\rm T} f(u, X^x_u)du\\
&+\int^t_{m-1} p(|X^x_u|^2+\delta)^\frac{p-2}{2}(X^x_u)^{\rm T} g(u, X^x_u)dB_u\\
&+ \int^t_{m-1} p(|X^x_u|^2+\delta)^\frac{p-2}{2}(X^x_u)^{\rm T} h^{ij}(u, X^x_u)\\
&+\frac{1}{2}\big(p(|X^x_u|^2+\delta)^\frac{p-2}{2} g^{ i}(u, X^x_u) g^{j}(u, X^x_u)\\[6pt] 
&+p(p-2)(|X^x_u|^2+\delta)^\frac{p-4}{2}(X^x_u)^{\mu}(X^x_u)^{\nu}g^{\mu i}(u, X^x_u)g^{\nu j}(u, X^x_u)\big)d\langle B^i B^j\rangle_u,\ q.s..
%
%&\leq \bar{\mathbb{E}}[|X^x_{n-1}|^p]+c_1\int^n_{n-1} \bar{\mathbb{E}}[|X^x_s|^p]ds+\bar{\mathbb{E}}\bigg[\int^n_{n-1} p|X^x_s|^{p-2}(X^x_s)^{\rm T} g(s, X^x_s)dB_s\bigg],
\end{align*}
From (\ref{cone}), we have
\begin{align}\label{er}
(|X^x_t|^2+\delta)^{\frac{p}{2}}\leq (|X^x_{m-1}|^2+\delta)^{\frac{p}{2}}&+\int^t_{m-1} pK(|X^x_u|^2+\delta)^\frac{p}{2}du+\int^t_{m-1} p(|X^x_u|^2+\delta)^\frac{p-2}{2}(X^x_u)^{\rm T} g(u, X^x_u)dB_u\notag\\
&+ c(\Sigma)\int^t_{m-1} pK(|X^x_u|^2+\delta)^\frac{p}{2}+\frac{pK}{2}\big(c(d)+c(n, d)|p-2|\big)(|X^x_u|^2+\delta)^\frac{p}{2}du\notag\\
&\leq (|X^x_{m-1}|^2+\delta)^{\frac{p}{2}}+c_1\int^t_{m-1} (|X^x_u|^2+\delta)^\frac{p}{2}du\notag\\
&+\int^t_{m-1} p(|X^x_u|^2+\delta)^\frac{p-2}{2}(X^x_u)^{\rm T} g(u, X^x_u)dB_u,
\end{align}
where $c_1>0$ is a constant  that depends only on $n$, $d$, $p$, $K$ and $\Sigma$. \\[6pt]
%Then, 
%\begin{align*}
%\bar{\mathbb{E}}[|X^x_n|^p]&=\bar{\mathbb{E}}\bigg[|X^x_{n-1}|^p+\int^n_{n-1} p|X^x_s|^{p-2}(X^x_s)^{\rm T} f(s, X^x_s)ds + \int^n_{n-1} p|X^x_s|^{p-2}(X^x_s)^{\rm T} h^{ij}(s, X^x_s) \\
%&+\frac{1}{2}\big(p|X^x_s|^{p-2} g^{ i}(s, X^x_s) g^{j}(s, X^x_s) +p(p-2)|X^x_s|^{p-4}(X^x_s)^{\mu i}(X^x_s)^{\nu j}g^{\nu j}(t, X^x_s)\big)d\langle B^i B^j\rangle_s\\
%&+\int^n_{n-1} p|X^x_s|^{p-2}(X^x_s)^{\rm T} g(s, X^x_s)dB_s\bigg]\\
%&\leq \bar{\mathbb{E}}[|X^x_{n-1}|^p]+c_1\int^n_{n-1} \bar{\mathbb{E}}[|X^x_s|^p]ds+\bar{\mathbb{E}}\bigg[\int^n_{n-1} p|X^x_s|^{p-2}(X^x_s)^{\rm T} g(s, X^x_s)dB_s\bigg],
%\end{align*}
%where the last inequality is in light of (\ref{cone}), Lemma \ref{lms}, the sub-additivity of $\bar{\mathbb{E}}[\cdot]$ and boundedness of the function $G$. 
\indent Then, we deduce by Lemma \ref{le5},
\begin{align*}
&\bar{\mathbb{E}}\bigg[\sup_{m-1\leq t\leq m}\bigg|\int^t_{m-1} p(|X^x_u|^2+\delta)^\frac{p-2}{2}(X^x_u)^{\rm T} g(u, X^x_u)dB_u\bigg|\bigg]\\
\leq\ &C\bar{\mathbb{E}}\bigg[\bigg(\int^m_{m-1} p^2(|X^x_u|^2+\delta)^{(p-2)}|(X^x_u)^{\rm T} g(t, X^x_u)|^2du\bigg)^{\frac{1}{2}}\bigg]\\
\leq\ &C\bar{\mathbb{E}}\bigg[\bigg(\big(\sup_{m-1\leq t\leq m}(|X^x_t|^2+\delta)^\frac{p}{2}\big)\int^m_{m-1} p^2 K(|X^x_u|^2+\delta)^\frac{p}{2}du\bigg)^{\frac{1}{2}}\bigg].
\end{align*}
On the other hand, by Young's inequality, we can always find a constant $c_2$ that depends on $c_1$ such that
\begin{align}\label{e58}
\notag\bar{\mathbb{E}}\bigg[\sup_{m-1\leq t\leq m}\int^t_{m-1} p(|X^x_u|^2&+\delta)^\frac{p-2}{2}(X^x_u)^{\rm T} g(u, X^x_u)dB_u\bigg]\\
&\leq \frac{1}{2}\bar{\mathbb{E}}\big[\sup_{m-1\leq t\leq m}(|X^x_t|^2+\delta)^\frac{p}{2}\big] + c_2p^2K\bar{\mathbb{E}}\bigg[\int^m_{m-1} (|X^x_t|^2+\delta)^\frac{p}{2}dt\bigg].
\end{align}
Thanks to the positivity of  $(|X^x|^2+\delta)^\frac{p}{2}$, we apply Fubini's theorem under each $\mathbb{P}\in\mathcal{P}_G$ to $\mathbb{E}^\mathbb{P}[\int^m_{m-1}(|X^x_t|^2+\delta)^\frac{p}{2}dt]$ and obtain 
\begin{equation}\label{e59}\bar{\mathbb{E}}\bigg[\int^m_{m-1} (|X^x_t|^2+\delta)^\frac{p}{2}dt\bigg]\leq \int^m_{m-1} \bar{\mathbb{E}}[(|X^x_t|^2+\delta)^\frac{p}{2}]dt\leq \max\{1, 2^{\frac{p}{2}-1}\}\int^m_{m-1} 
(\bar{\mathbb{E}}[|X^x_t|^p]+\delta^\frac{p}{2})dt.
\end{equation}
%although $(|X^x|^2+\delta)^\frac{p}{2}$ may be not a $M^1_*([0, T])$ element.
We conclude from (\ref{stacond})-(\ref{e59}) that
\begin{align*}
\bar{\mathbb{E}}[\sup_{m-1\leq t\leq m}|X^x_t|^p]&\leq \bar{\mathbb{E}}[\sup_{m-1\leq t\leq m}(|X^x_t|^2+\delta)^{\frac{p}{2}}]\\
&\leq 2\max\{1, 2^{\frac{p}{2}-1}\}\bar{\mathbb{E}}[|X^x_{m-1}|^p]+c_3\int^m_{m-1}\bar{\mathbb{E}}[|X^x_t|^p]dt+c'_3\delta^\frac{p}{2}
\end{align*}
where $c_3$, $c'_3>0$ depends only on $c_1$, $c_2$ and $p$. Letting $\delta\rightarrow 0$, there exists a constant $c_4$ depends only on $c_3$, $c'_3$ and $x$, such that 
$$
\bar{\mathbb{E}}[\sup_{m-1\leq t\leq m}|X^x_t|^p]\leq c_4\exp(-\lambda(m-1)).
$$
\indent Applying Markov's inequality (Lemma \ref{MI}) gives, for arbitrary $\varepsilon \in (0, \lambda)$, 
$$
\bar{C}\bigg\{\sup_{m-1\leq t\leq m}|X^x_t|^p>e^{-(\lambda -\varepsilon )(m-1)}\bigg\}\leq c_4 \exp(-\varepsilon (m-1)),
$$
which implies
$$
\sum^\infty_{n=1}\bar{C}\bigg\{\sup_{m-1\leq t\leq m}|X^x_t|^p>e^{-(\lambda -\varepsilon )(m-1)}\bigg\}<\infty.
$$
From Borel-Cantelli lemma (Lemma \ref{BCL}), we see that for q.s. every $\omega$, there exists $N_0:=N_0(\omega)$, such that for $m\geq N_0(\omega)$ and $t\in [m-1, m]$,
$$
\frac{1}{t}\log(|X^x_t|)=\frac{1}{pt}\log(|X^x_t|^p)< -\frac{(\lambda-\varepsilon)(m-1)}{pm},
$$
Passing to limit as $\varepsilon \rightarrow 0$ and $m\rightarrow \infty$, we obtain the desired result. \hfill$\square$\\[6pt]
\indent The next property of the paths of the $G$-Brownian motion follows from the theorem above.
%his theorem, we can deduce a corollary that describes the pathwisely asymptotic property of $\frac{B_t}{t}$. 
\begin{corollary}
$$
\frac{|B_t|}{t}\rightarrow 0,\ {\rm{as}}\ t\rightarrow \infty,\ q.s..
$$
\end{corollary}
\noindent {\bf Proof.} Without loss of generality, we only prove this property for scalar $G$-Brownian motion. Consider a particular case of Example \ref{pg}, i.e. in (\ref{leqe}), $\alpha=-1$, $\beta=\frac{1}{2}$ and $\gamma =1$. Then, we calculate for some $0<p<\frac{2}{\overline{\sigma}^2}$,
$$
\bar{\mathbb{E}}[|X^x_t|^p]\leq |x|^p\exp\big(-p(1-\frac{p\overline{\sigma}^2}{2})t\big).
$$
By Theorem \ref{mows}, we obtain
$$
\limsup_{t\rightarrow \infty}\frac{1}{t}\log(|X^x_t|)= \limsup_{t\rightarrow \infty}\frac{1}{t}\log(|\exp(-t+B_t)|)=-1+\limsup_{t\rightarrow \infty} \frac{B_t}{t}\leq-{(1-\frac{p\overline{\sigma}^2}{2})},
$$
which implies $\limsup_{t\rightarrow \infty} \frac{B_t}{t}\leq \frac{p\overline{\sigma}^2}{2}$. On the other hand, we consider $\alpha=-1$, $\beta=\frac{1}{2}$ and $\gamma =-1$ in (\ref{leqe}) and obtain $\liminf_{t\rightarrow \infty} \frac{B_t}{t}\geq -\frac{p\overline{\sigma}^2}{2}$. Letting $p\rightarrow 0$, we complete the proof. \hfill$\square$
\indent Now we turn to discuss the sufficient condition for the moment instability. In order to ensure the solvability of (\ref{gsde1}) in the following theorem, we always assume Lipschitz conditions on its coefficients.
\begin{theorem}
Consider the $G$-stochastic system (\ref{gsde1}), whose coefficients are deterministic functions continuous in $t$ and Lipschitz in $x$. If there exists a function $V(t, x)\in\mathcal{C}^{1, 2}([0, \infty)\times \mathbb{R}^n)$ such that for all $t\geq 0$, some positive constant $c_1$, $c_2$, and some $q>0$, 
\begin{equation}\label{50k}
c_1 |x|^{q}\leq V(t, x)\leq c_2|x|^{q}
\end{equation}
and 
\begin{equation}\label{03}
%\underline{\mathcal{L}}(V)(t, x):=-\mathcal{L}(-V)(t, x)\geq \lambda V(t, x).
%\lambda V(t,x ). 
\mathcal{L}V(t, x)\geq \lambda V(t, x).
\end{equation}
Then, the trivial solution is exponentially $q$-instable, namely,
\begin{equation*}
\bar{\mathbb{E}}[|X^{s, x}_t|]\geq C|x|^{q}\exp(\lambda(t-s)).
\end{equation*}
%%of the $G$-stochastic system (\ref{gsde1})
%\item the trivial solution is exponentially $p$-stable, if there exists a $\lambda>0$ such that for all $(t, x)\in ([0, \infty)\times \mathbb{R}^n)$,
%\end{enumerate}
\end{theorem}
\noindent{\bf Proof:} Without loss of generality, we prove only the case $s=0$. In the light of the property of the solution of a Lipschitz GSDE, one can define a sequence of stopping times by $
\tau_N:=\inf\{t: X^x_t\geq N\}, 
$
such that $X^{x}{\bf 1}_{[0, \tau_N]}\in M^p_*([0, T]; \mathbb{R}^n)$, for any $p\geq 1$.
 Applying $G$-It\^o's formula to $\Phi(t\wedge \tau_N, X_{t\wedge\tau_N}):= -\exp(-\lambda(t\wedge \tau_N))V(t\wedge\tau_N, X_{t\wedge\tau_N})$, one can follow the proof of Theorem \ref{lf} and obtain
$$
\bar{\mathbb{E}}[\Phi(t, X^x_{t})]\leq \liminf_{N\rightarrow \infty}\bar{\mathbb{E}}[\Phi(t\wedge \tau_N, X^{x}_{t\wedge \tau_N})]\leq \Phi(0, x).
$$
Thus, 
$$
\bar{\mathbb{E}}[V(t, X^x_{t})]\geq V(0, x)\exp(\lambda t),
$$
which implies 
$$
\bar{\mathbb{E}}[|X^x_{t}|^q]\geq \frac{c_1}{c_2}|X^x_{t}|^q\exp(\lambda t).
$$
\hfill$\square$
\begin{example}
Consider the linear GSDE (\ref{ex}), a sufficient condition to ensure that the trivial solution is $q$-unstable can be given by the following system of LMIs: for some $\alpha\in\mathbb{R}$,
\begin{align*}
\left\{\begin{array}{l}
2PH+C^{\rm T}PC\geq \alpha I_n;\\[3pt]
2PF-I_n\geq -G(\alpha) I_n.
\end{array}\right.
\end{align*}
\end{example}
%necessary condition for the moment stability of the trivial solution. In order to construct a function $V$ that is sufficient smooth, we will adopt the following condition:\\[6pt]
%{\bf (H5')} The coefficients $f(\cdot, \cdot)$, $h^{ij}(\cdot, \cdot)$, $g^{j}(\cdot, \cdot): [0, T]\times\mathbb{R}^n\rightarrow \mathbb{R}^n$ are deterministic functions continuous in $t$ and have continuous bounded derivatives with respect to $x$ up to second order. 
%\begin{theorem}
%Consider the $G$-stochastic system (\ref{gsde1}) with (H4') and (H5'). If the trivial solution is exponentially $p$-stable, namely,
%for all $x\in \mathbb{R}^n$,
%\begin{equation}\label{stacond1}
%\bar{\mathbb{E}}[|X^{s, x}_t|^p]\leq |x|^p \exp(-\lambda (t-s)),
%\end{equation}
%where $\lambda>0$. Then, there exists a function $V\in \mathcal{C}^{1, 2}([0, T]\times\mathbb{R}^n)$ satisfying (\ref{54k}), (\ref{3}) and 
%\begin{equation}
%\bigg|\frac{\partial V}{\partial x^\mu}\bigg|<c_1|x|^{p-1},\ \bigg|\frac{\partial^2 V}{\partial x^\mu\partial x^\nu}\bigg|<c_2|x|^{p-2}.
%\end{equation}
%\noindent If the trivial solution is exponential $p$-stable, for some $p >0$, namely, for all $x\in \mathbb{R}^n$,
%\begin{equation}\label{stacond}
%\bar{\mathbb{E}}[|X^x_t|^p]\leq |x|^p \exp(-\lambda t),
%\end{theorem}%locally Lipschitz in $x$,
%\noindent {\bf Proof:} First, the assumption (H4'') implies the uniformly Lipschitz condition (H1) and since all coefficients are deterministic, the solution $X^x$ of (\ref{gsde1}) stays in $M^p_G([0, T]; \mathbb{R}^n)\subset M^p_*([0, T]; \mathbb{R}^n)$, for any $p\geq 2$. 
\section{Appendix}
\noindent In the Appendix, we study the property of stochastic integrals with respect to the quadratic variation process $\langle B, B\rangle$, and we obtain a generalized result of Proposition 1.4 in \S IV-1 of Peng \cite{P3}, for a $\eta\in M_w^1([0, T], \mathbb{S}^d)$. 
%\begin{proposition}
%There exists a $\mathcal{B}([0, t])\times \mathcal{F}_t$-progressively measurable process $\hat{a}$ taking values in $\Sigma$, such that for any $\eta=(\eta^{ij})^{d}_{i, j=1}\in M^1_*([0, T]; \mathbb{S}^d)$ and $t\in [0, T]$,
%$$
%\int^t_0 \eta^{ij}_s d\langle B^i B^j\rangle_s=\int ^t_0 {\rm tr} (\hat{a}_s\eta_s) ds,\ \mathbb{P}-a.s.\ {\rm for\ each}\ \mathbb{P}\in \mathcal{P}_1.
%$$
%\end{proposition}
%\noindent{\bf Proof:}\\
%\hfill$\square$\\[6pt]
%\noindent The following lemma can be regarded as a direct corollary of the proposition above:
\begin{lemma}
Let $\eta=(\eta^{ij})^{d}_{i, j=1}\in M^1_*([0, T]; \mathbb{S}^d)$ and 
$
M_t:=\int^t_0 \eta^{ij}_s d\langle B^i B^j\rangle_s-\int^t_0 2G(\eta_s) ds$. Then, for each $t\in [0, T]$, $M_t\leq 0$, q.s., thus $\bar{\mathbb{E}}[M_t]\leq 0$. 
\end{lemma}
\noindent {\bf Proof:} Without loss of generality, we prove only the case $d=1$. By the definition, we know that $\langle B\rangle$ has continuous and increasing paths on $[0, T]$, q.s.. From Corollary 5.4 in \S III-5 of Peng \cite{P3}, we have, for any $0\leq t_1\leq t_2\leq T$, in the sense of $L^1_G(\Omega_T)$ (consequently in the sense of q.s.),
\begin{equation}\label{zj}
\underline{\sigma}^2(t_2-t_1)\leq \langle B\rangle_{t_2}-\langle B\rangle_{t_1}\leq \overline{\sigma}^2(t_2-t_1).
\end{equation}
Thus, one can pick up a $\bar{C}$-null set $A$, such that for each $\omega\in A^c$, for each $t_1$, $t_2\in \mathbb{Q}\cap [0, T]$, (\ref{zj}) holds true, furthermore, by the continuity of the path, (\ref{zj}) indeed holds true for all  $t_1$, $t_2\in [0, T]$. Therefore, $\langle B\rangle_\cdot(\omega)$ induce a finite measure on $[0, T]$, which is absolutely continuous w.r.t. the Lebesgue measure. We denote by $\hat{a}_\cdot(\omega)$ the R-N derivative $\frac{d\langle B\rangle_t(\omega)}{dt}$, which is in $[\underline{\sigma}^2, \overline{\sigma}^2]$, $\lambda$-a.s.. As a result,
$$M_t(\omega):=\int^t_0 \eta_s(\omega) d\langle B\rangle_s(\omega)-\int^t_0 2G(\eta_s(\omega)) ds\leq 0,\ q.s..$$
Taking $\bar{\mathbb{E}}[\cdot]$ on both sides, we obtain the desired result.\hfill$\square$
\begin{remark}
In \S 2.4.1 of Lin \cite{L2}, the following fact is proved: for $\eta\in M^1_*([0, T])$, the stochastic integral of $\eta$ w.r.t. $\langle B\rangle$ on $[0, T]$ can be defined in both the Lebesgue-Stieltjes sense and the sense of a contracting mapping $M^1_*([0, T])\rightarrow L^1_*([0, T])$; these two definitions are indistinguishable in the $M^1_*(\Omega_T)$ sense.
\end{remark}
\indent Furthermore, we can generalize this lemma to the case that $\eta=(\eta^{ij})^{d}_{i, j=1}\in M^1_w([0, T]; \mathbb{S}^d)$.
\begin{lemma}\label{lms}
Let $\eta=(\eta^{ij})^{d}_{i, j=1}\in M^1_w([0, T]; \mathbb{S}^d)$ and 
$
M_t:=\int^t_0 \eta^{ij}_s d\langle B^i B^j\rangle_s-\int^t_0 2G(\eta_s) ds$. Then, for each $t\in [0, T]$, $\bar{\mathbb{E}}[M_t]\leq 0$. 
\end{lemma}
\noindent {\bf Proof:}
The definition of $M$ is well posed because of (\ref{haha}) and the Lipschtiz continuity of the function $G$. Suppose $\eta$ is associated with the sequence of stopping time $\{\sigma_m\}_{m\in \mathbb{N}}$. Then, 
\begin{align*}
M_{t\wedge\sigma_m\wedge T}&=\int^{t \wedge \sigma_m\wedge T}_0 \eta^{ij}_s d\langle B^i B^j\rangle_s-\int^{t\wedge \sigma_m\wedge T}_0 2G(\eta^{ij}_s) ds\\
&=
\int^{t}_0 \eta^{ij}_s {\bf 1}_{[0, \sigma_m\wedge T]}(s)d\langle B^i B^j \rangle_s-\int^{t}_0 2G(\eta_s{\bf 1}_{[0, \sigma_m\wedge T]}(s)) ds\leq 0,\ q.s..
\end{align*}
Fatou's Lemma implies
$$
\bar{\mathbb{E}}[M_t]=\bar{\mathbb{E}}[\lim_{m\rightarrow \infty}M_{t\wedge\sigma_m\wedge T}]\leq \liminf_{m\rightarrow \infty}\bar{\mathbb{E}}[M_{t\wedge\sigma_m\wedge T}]\leq 0.
$$
\hfill$\square$\\[6pt]
\noindent\textbf{Acknowledgement}\ The authors expresses special thanks to Prof. Ying HU for his useful suggestions and also to Prof. Shige PENG for his constructive advice to improve this paper. 

\end{document}